%Paper: q-alg/9506019
%From: Yao-zhong Zhang <yzzhang@yukawa.kyoto-u.ac.jp>
%Date: Fri, 23 Jun 1995 15:41:34 +0900
%Date (revised): Wed, 26 Jul 1995 16:36:31 +0900

\input amstex
\magnification=\magstep1
\documentstyle{amsppt}
\NoBlackBoxes
\font\BIG=cmbx12 scaled\magstep1

\define\vs{\vskip.2cm}
\define\ve{\varepsilon}
\define\La{\Lambda}
\define\la{\lambda}

\define\ga{\gamma}
\define\ot{\otimes}
\define\Z{\Bbb Z}
\define\om{\omega}
\define\bal{\beta}
\define\al{\alpha}

\define\tr{\Delta}
\pageno=1
\document
\baselineskip=6mm plus 2pt
\vskip.2cm
\TagsOnRight
%\pagewidth { 130 mm }
%\pageheight {180 mm}
%\hcorrection {6 mm}
\parindent 8 mm

\rightline{Lett. Math. Phys. in press}
\rightline{YITP/K-1112}
\rightline{
q-alg/9506019
}

\vskip 1cm
\centerline{\BIG `` EIGENVALUES OF CASIMIR INVARIANTS FOR}
\centerline {\BIG TYPE I QUANTUM SUPERALGEBRAS ''}
\vskip 1cm
\centerline {\bf by}
\vskip 1cm
\centerline {\bf Mark D. Gould, Jon R. Links}
\+\cr
\centerline {Department of Mathematics,
            University of Queensland,
            Queensland, 4072, Australia}
\+\cr
\centerline {and}
\+\cr
\centerline {\bf Yao-Zhong Zhang}
\+\cr
\centerline {Yukawa Institute for Theoretical Physics,
             Kyoto University,
              Kyoto 606, Japan}
\+\cr
\centerline {March, 1995 }
\+\cr
\+\cr
\vskip.3cm
\noindent
We present the eigenvalues of  the Casimir
invariants for  the type I quantum superalgebras on
any irreducible highest weight module.
\+\cr
\+\cr
\centerline {Mathematics Subject Classifications 81R10, 17B37, 16W30}
\+\cr
\+\cr

\vfil
\eject
\heading 1. Introduction  \endheading
\+\cr

Quantum algebras \cite {1} are well known for their role in solving the
Yang-Baxter equation and representing the symmetries of the associated
solvable models. Likewise, their $\Bbb Z_2$-graded counterparts quantum
superalgebras \cite {2, 3, 4} play a similar role in relation to supersymmetric
solvable models. There has been significant interest in the study of
these models, particularly those describing systems of correlated
electrons which are of importance in condensed matter physics. Such
examples are the supersymmetric t-j model \cite {5} and its quantum
analogue \cite {6}, the supersymmetric extended Hubbard model \cite {7}
and the model of Bracken et. al. \cite {8}. All of these models are
exactly solvable in one dimension.

In the investigation of any model possessing a quantum (super)algebra
as a symmetry, it is desirable to have a well developed representation
theory. An important part in this pursuit is to study the central
elements or Casimir invariants. There are established techniques
for the construction of Casimir invariants, both for quantum algebras
\cite {9} and quantum superalgebras \cite {10, 11}. Here we follow the
construction of \cite {10} where the Casimir invariants are obtained
from tensor operators derived from the universal $R$-matrix and an {\it
arbitrary} reference representation. Such general invariants are
important in the computation of link invariants in knot theory \cite
{12}. However a formula for the eigenvalues of these invariants when
acting on finite dimensional irreducible representations has never been
derived. In this Letter we will present such a formula  for the type I
quantum superalgebras consisting of $U_q(gl(m|n))$ and $U_q(osp(2|2n))$.
The proof of the formula is both lengthy and detailed and will be
deferred to a separate publication. Here we will merely report our
results which are new even in the classical case ($q\rightarrow 1$
limit).  Our results are in complete agreement  with studies of some
particular cases \cite {13, 14},  but  are more general and widely
applicable.

Although our results apply only for the type I quantum superalgebras,
they are the most interesting for the following reasons. Firstly, they
admit finite dimensional unitary representations \cite {15} which have
physical importance where unitarity is a requirement (e.g. see \cite
{8}). Secondly, the type
I quantum superalgebras admit one parameter families of representations
which have interesting applications such as providing solutions to the
Yang-Baxter equation with extra spectral parameters, though in a
non-additive form \cite {16}. Another application is in the construction
of two variable link invariants as indicated in \cite {17}. However, there
have always been technical difficulties in evaluating such link
invariants. The results reported in  the present Letter now permit a unified
construction and evaluation of these invariants \cite {18}.

\heading 2. Preliminaries\endheading
\vs

Let $g$ denote a basic classical
Lie superalgebra of rank $l+1$ with the usual generators
$\{e_i,\,f_i,\,h_i\}_{i=0}^l$. Let $\{\alpha_i\}_{i=0}^l$ be the
distinguished set of simple
roots of $g$ in the sense of Kac \cite {19} and let $(\,\,\, ,\,\,\,)$ be a
fixed invariant bilinear
form on $H^*$, the dual of the Cartan subalgebra $H$ of $g$. We also let
$\Phi^+=\Phi^+_0\cup \Phi^+_1$ denote the full set of positive roots with
$\Phi_0^+$ (resp. $\Phi_1^+$) the subset of even (resp. odd) positive roots.
Throughout,
we adopt the convention that $\alpha_0$ denotes the unique odd simple
root. Associated
with $g$ one can define the quantum superalgebra $U_q(g)$
($q$ is assumed  not   a root of unity) which has the
structure of a $\Z_2\snug$-graded quasi-triangular Hopf algebra
\cite{3}. We will not give the full defining relations of $U_q(g)$ here
and refer to  \cite{4} for details. We note however that
$U_q(g)$ has a co-product structure given by
$$\tr(q^{\pm\frac12 h_i})=q^{\pm\frac12 h_i}\ot q^{\pm\frac12
h_i},\quad \tr(x)=x\ot q^{-\frac12 h_i}+q^{\frac12 h_i}\ot x,\quad
x=e_i,\,f_i$$
which is extended to an algebra homomorphism to all of $U_q(g)$ in the
usual way. It is important to point out that the multiplication rule for
the tensor product is defined for homogeneous elements $a,b,c,d\in
U_q(g)$ by
$$(a\ot b)(c\ot d)=(-1)^{[b][c]}(ac\ot bd)\tag1$$
and extended linearly to all of $U_q(g)\ot U_q(g)$. Here $[a]\in\Z_2$
denotes the degree of the homogeneous element $a\in U_q(g)$, which is
defined for the elementary generators by
$$[h_i]=0,\quad [e_i]=[f_i]\equiv [i]=\delta_{i0},\qquad\forall\, 0\leq
i\leq l, $$
and extended to all homogeneous elements of $U_q(g)$ through
$$[ab]=[a]+[b](mod\, 2),\qquad\forall\,a,b\in U_q(g).$$

\vs
\define\ra{\rightarrow}
The twist map $T\:U_q(g)\ot U_q(g)\ra U_q(g)\ot U_q(g)$ is defined by
$$T(a\ot b)=(-1)^{[a][b]}b\ot a\tag 2$$
\define\otr{\overline{\Delta}}
for all homogeneous $a,b\in U_q(g)$: we set $\otr=T.\tr$. Then there
exists a canonical element $R\in U_q(g)\ot U_q(g)$ called the universal
$R\snug$-matrix which is even and invertible and satisfies the following
relations \cite{3}
$$\gather
R\tr(a)=\otr(a)R,\qquad\forall a\in U_q(g), \tag 3\\
(\tr\ot I)R=R_{13}R_{23},\quad(I\ot\tr)R=R_{13}R_{12}, \tag 4\endgather
$$
where we have adopted the conventional notation. From equations (3) and
(4) it follows that the universal $R\snug$-matrix satisfies the graded
Yang-Baxter equation
$$R_{12}R_{13}R_{23}=R_{23}R_{13}R_{12}.\tag 5$$
We emphasize that multiplication of the tensor products are to obey
equation (1).

Let $\rho\in H^*$ denote the graded half sum of positive roots of $g$
and let $h_{\rho}$ denote the unique element of $H$ defined by $\alpha_i
(h_{\rho})=(\rho,\,\alpha_i),\quad\forall \alpha_i\in H^*$. We recall
from  \cite {10} the following result.

\proclaim{Theorem 1} Let $\pi$ be a fixed, but arbitrary, finite
dimensional representation of $U_q(g)$ with representation space $V$ and
set
$$\Delta_{\pi}=(\pi \otimes I)\Delta.$$
If $w\in  \,End\, V\otimes U_q(g)$ satisfies
$$\Delta_{\pi}(a)w=w\Delta_{\pi}(a),\qquad\quad\forall a \in U_q(g) \tag
{6} $$
then
$$s\tau_q(w)=(str \otimes I)(\pi(q^{2h_{\rho}})\otimes I)w $$
belongs to the centre of $U_q(g)$, where $str$ denotes the supertrace.
\endproclaim

Theorem 1 enables a family of Casimir invariants to be constructed for
$U_q(g)$ for any reference module $V$ utilizing the universal $R$-matrix.
Defining $R^T=TR$, it is
clear from (3) that
$$R^TR\Delta (a)=\Delta (a)R^TR,\qquad\quad \forall a \in U_q(g).$$
Setting
$$A=(q-q^{-1})^{-1}(\pi\otimes I)(I\otimes I-R^TR)$$
then $A^l,\, l\in\Bbb Z^+$, satisfies (6). We thus obtain the family of
Casimir invariants
$$C_l=s\tau_q (A^l).   \tag 7 $$
In the limit $q\rightarrow 1$ these give rise to a family of Casimir
invariants for the classical Lie superalgebra $g$.

The preceeding discussion applies for any
quantum superalgebra, not only those of type I.
Let $V(\mu )$ denote a  finite dimensional irreducible $U_q(g)$ module
of highest weight $\mu\in D^+$ where $D^+\subset H^*$ is the set of
dominant weights. For type I quantum superalgebras, $\mu\in D^+$
if and only if \cite {19, 20}
$$<\mu,\,\al_i>=\frac {2(\mu,\,\al_i)}{(\al_i,\,\al_i)}\in \Bbb Z^+,
\qquad 1\leq i\leq l,$$
while $(\mu,\,\al_0)$ can take arbitrary complex values.
When acting on $V(\mu)$ the invariants $C_l$ act as
scalar multiples of the identity operator (Schur's lemma), which we
denote by $\chi_{\mu}(C_l)$. A general formula for these eigenvalues is
unknown. In the next section we will present such a formula for the type
I quantum superalgebras.

\heading 3. Eigenvalue Formula \endheading

Hereafter $U_q(g)$ is assumed to be of type I.  Recall that for type I
quantum superalgebras we say that $\mu\in H^*$ is {\it typical} if \cite
{19, 20}
$$(\mu +\rho,\,\alpha)\neq 0,\qquad \forall\alpha\in\Phi^+_1,$$
and {\it atypical} otherwise.  We define the following.
\proclaim{Definition 1} Let $\Pi(\La)$ denote the weight spectrum of
$V(\La)$. We say that $\La$ is {\it subordinate} to $\mu\in D^+$ if $\,\forall
\la\in \Pi(\La),\,\,\mu
+\la$ is dominant. Moreover we say that $\La$ is {\it typically
subordinate} to $\mu $ if $\mu,\,\mu +\la$ are all typical and dominant.
We denote the set of such $\mu$ by
$D^+_{\La}$.
\endproclaim

Let ${\la_i}$ denote the distinct weights in
the reference module $V(\La)$ with $[\la_i]$ the
degree of $\la_i$.  For $\mu\in
D^+_{\La}$ we have the tensor product decomposition
$$V(\La)\ot V(\mu)=\bigoplus_i m_i V(\mu +\la_i), \tag 8$$
where $m_i$ denotes the multiplicity of $V(\mu +\la_i)$. The above
decomposition is necessarily completely reducible since any typical
module splits in a finite dimensional representation \cite {19}. Also
$m_i$ is the same as the multiplicity of $\la_i$ occurring in $\Pi(\La)$
 (see lemma A2 in the appendix).
If
$$P[i]\equiv P[\mu +\la_i] $$
denote central projections onto the isotypic components
$\overline{V(\mu +\la_i)}\equiv m_iV(\mu+\la_i)$ we have the spectral
decomposition \cite {10}
$$A^l=\sum_i  \left [\bal_i(\mu)\right ]^l P[i],$$
where the roots $\bal_i(\mu)$ are given by
$$\bal_i(\mu)=\frac {1-q^{-(\la_i,\la_i+2\mu +2\rho )+(\La,\La +2\rho)}}
{q-q^{-1}} .$$

Since the $P[i]$ are central projections, we may invoke theorem 1 to
construct the central elements
$$\ga[i]=s\tau_q(P[i]),$$
which leads to
$$\chi_{\mu}(C_l)=\sum_i\left  [\bal_i(\mu)\right ]^l\chi_{\mu}\left
(\ga[i]\right ).$$
Recently, we have determined the quantities $\chi_{\mu}\left (\ga[i]
\right )$ \cite {21}. The details of the derivation of our results are
beyond the scope of this Letter. Here we merely wish to present our
results.

\proclaim{Proposition 1} For any $\overline{V(\mu +\la_i)}\subset
V(\La)\otimes V(\mu),\,\la_i\in\Pi(\La)$ such that $\mu,\,\mu +\la_i$ are
both typical and dominant then
$$\chi_{\mu}\left (\ga[i]\right )= (-1)^{[\la_i]}m_i\prod_{\al\in \Phi^+_0}
\frac {[(\mu+\la_i+\rho,\,\al)]_q}{[(\mu+\rho,\,\al)]_q}\prod_{\al\in
\Phi_1^+} \frac {[(\mu+\rho,\,\al)]_q}{[(\mu+\la_i+\rho,\,\al)]_q},
$$
where
$$[x]_q=\frac{q^x-q^{-x}}{q-q^{-1}}.$$
\endproclaim

With the aid of proposition 1, we can determine the eigenvalues of the
Casimir invariants $C_l$ on all irreducible finite dimensional modules
$V(\mu)$ such that $\La$ is typically subordinate to $\mu$. In view of the
following, we can extend our formula to all modules $V(\mu)$.

\proclaim{Proposition 2}
Let $\La\in D^+ $ be fixed but arbitrary. If $f$ is a polynomial function
on $H^*$  satisfying  $f(\nu)=0, \quad\forall \nu\in D^+_\La, $ then $f$
vanishes identically.    \endproclaim
\noindent
The proof of proposition 2 is left to the appendix.

\parindent  8mm
Consider the quantity
$$f(\mu) =\chi_{\mu}(C_l)-\sum_i\left [\bal_i(\mu)\right
]^l\chi_{\mu}\left (\ga[i]\right ),$$
with $\chi_{\mu}\left (\ga [i]\right )$ as in proposition 1. It is apparent
that
$$f(\mu)=0,\qquad\forall \mu\in D^+_{\La}. $$
Expanding $f(\mu)$ into a power series in $\eta =ln\,q$ we obtain
$$f(\mu)=\sum_k f_k(\mu)\eta^k,$$
where $f_k(\mu )$ is a polynomial function on $H^*$. It follows that
$$f_k(\mu)=0,\qquad\forall \,\mu\in D^+_{\La}.$$
In view of proposition 2 we have
$$f_k(\mu)=0,\qquad\forall\,\mu\in H^*, $$
from which we deduce the following.
\proclaim{Proposition 3} The eigenvalues of the Casimir invariants (7)
acting on the irreducible module $V(\mu),\quad \forall\mu\in D^+,$ are given by
$$\chi_{\mu}(C_l)=\sum_i(-1)^{[\la_i]}m_i\left [\bal_i(\mu)\right ]^l
\prod_{\al\in\Phi^+_0}\frac{[(\mu+\la_i+\rho,\,\al)]_q}{[(\mu+\rho,
\,\al)]_q}\prod_{\al\in\Phi^+_1}\frac{[(\mu+\rho,\,\al)]_q}{[(\mu+\la_i
+\rho,\,\al)]_q}.  \tag 9$$
\endproclaim

The above eigenvalue formula in fact holds on any $U_q(g)$ module
admitting an infinitesimal character $\chi_{\mu},\,\mu\in H^*$. However
for atypical $\mu +\la_i$, proposition 3 as presented is undefined.
Nevertheless, it is important to observe that $\chi_{\mu}(C_l)$ is in fact a
polynomial function in $q^{\pm(\mu,\,\al_i)},\,i=0,1,....,l.$ Hence in
principle eq. (9) may be expanded to yield a well defined expression.

In the limit $q\rightarrow 1$ we obtain the following result for the
classical Lie superalgebra case:
$$\chi_{\mu}(C_l)=\sum_i(-1)^{[\la_i]}\left [\bal_i(\mu)\right ]^l\prod_{\al
\in\Phi_0^+}\frac{(\mu +\la_i+\rho,\,\al)}{(\mu +\rho,\,\al)}\prod_{\al\in
\Phi_1^+}\frac{(\mu +\rho,\,\al)}{(\mu+\la_i+\rho,\,\al)},$$
where
$$\bal_i(\mu)=\frac 12(\la_i,\,\la_i+2\mu+2\rho)-\frac12(\La,\,\La
+2\rho).$$
This result is also new. It generalizes to arbitrary
reference representations the results of \cite {14}.

\+\cr
\+\cr
\heading 4. Examples \endheading
\+\cr
We will now illustrate our results for the case when $\La$ is the vector
representation. Let us first consider $U_q(gl(m|n))$. We choose $\{\ve_i
\}_{i=1}^{m+n}$ as a basis for $H^*$ with the $\Bbb Z_2$-gradation
$$[\ve_i]\equiv [i]=\cases 0, & 1\leq i\leq m, \\
			   1, & m<i\leq m+n,     \endcases
$$ and equipped with the invariant bilinear form
$$(\ve_i,\,\ve_j)=(-1)^{[i]}\delta_{ij}.$$
The sets of even and odd positive roots are respectively given by
$$\align
\Phi^+_0&=\{\ve_i-\ve_j|1\leq i<j\leq m+n,\,[i]=[j]\}, \\
\Phi^+_1&=\{\ve_i-\ve_j|1\leq i<j\leq m+n,\,[i]\neq [j]\},
\endalign $$
in terms of which  the graded half sum of positive roots is expressed as
$$\rho=\frac 12\sum_{i=1}^m(m-n-2i+1)\ve_i-\frac 12\sum_{j=1}^{n}
(m+n-2j+1)\ve_{m+j}.$$

The vector representation of $U_q(gl(m|n))$ has highest weight
$\La=\ve_1$ and the full weight spectrum is
$\Pi(\ve_1)=\{\ve_i\}_{i=1}^{m+n}$. Applying formula (9) for the
eigenvalues of the Casimir invariants (7) (with $\La=\ve_1$)
acting on the module $V(\mu)$  yields
$$
\chi_{\mu}(C_l)=\sum_{i=1}^{m+n} (-1)^{[i]}\left [\bal_i(\mu)\right ]^l
\prod_{j\neq i}^{m+n}\frac {q^{-(\ve_j,\ve_j)}\bal_i(\mu)-q^{(\ve_j,\ve_j
)}\bal_j(\mu) + (-1)^{[j]}   }
{\bal_i(\mu)-\bal_j(\mu)}   $$
where the roots $\bal_i(\mu)$ are given by
$$
\bal_i(\mu)=\frac {1-q^{-(\ve_i,\ve_i+2\mu +2\rho)+m-n}}{q-q^{-1}}.$$
The above eigenvalue formula is in agreement with the results of \cite
{13} obtained by different means. In the $q\rightarrow 1$ limit the
above formula coincides with those of \cite {14} taking into account the
different conventions adopted.

Next we consider the case $U_q(osp(2|2n))$. We choose the basis
$\{\ve_i\}_{i=0}^n$ for $H^*$ with the invariant bilinear form
$$(\ve_i,\,\ve_j)=(-1)^{[i]}\delta_{ij},$$
where $\ve_0$ is the unique odd basis element. The sets of even and odd
positive roots are respectively given by
$$\align
\Phi^+_0&=\{\ve_i\pm\ve_j|1\leq i<j\leq n\}\cup \{2\ve_i|1\leq i\leq
n\},   \\
\Phi^+_1&=\{\ve_0\pm \ve_i|1\leq i\leq n \}, \endalign  $$
while the graded half sum of positive roots reads
$$\rho =-n\ve_0+\sum_{i=1}^n(n-i+1)\ve_i.$$

The vector representation of $U_q(osp(2|2n))$ has highest weight
$\La=\ve_0$ and weight spectrum $\Pi (\ve_0)=\{\pm\ve_i\}_{i=0}^n$.
Applying the eigenvalue formula (9) yields
$$\align
\chi_{\mu}(C_l)=&\sum_{i=-n}^n (-1)^{[i]}\left [\bal_i(\mu)\right ]^l
\frac{q^{-1-(\ve_i,\ve_i)}\bal_i(\mu)-q^{1+(\ve_i,\ve_i)} \bal_{-i}
(\mu)+q+(-1)^{[i]}q^{-(\ve_i,\ve_i)}}{\bal_i(\mu)-\bal_{-i}(\mu)} \\
&\times \prod_{j\neq
\pm i} \frac {q^{-(\ve_j,\ve_j)}\bal_i(\mu)-q^{(\ve_j,\ve_j)}\bal_j(\mu)
+(-1)^{[j]} }{\bal_i(\mu)-\bal_j(\mu) }
\endalign$$
where
$$
\bal_i(\mu)= \frac {1-q^{-(\ve_i,\ve_i+2\mu +2\rho) +2n-1}}{q-q^{-1}}.
$$
In the  above formula it should be understood that
$$\ve_{-i}\equiv-\ve_i,\qquad [-i]\equiv [i].$$

For generic values of $q$ the above eigenvalue formula is a new result.
As $q\rightarrow 1$ it agrees with the formulae presented in \cite {14}.

\heading 5. Conclusion \endheading

In this Letter we have presented a general eigenvalue formula for the
Casimir invariants of the type I quantum superalgebras when acting
on  irreducible finite dimensional modules. Our formula applies
to arbitrary reference modules and the eigenvalues are expressed in
terms of the highest weight labels. The construction of these
Casimir invariants is given in an earlier publication \cite {10}. Our
results are new even in the classical $q\rightarrow 1$ limit and will
have application in a variety of areas including the evaluation of link
invariants \cite {18}. As examples we have explicitly evaluated these
eigenvalues in the case when the reference module is the vector module.
For the case of $U_q(gl(m|n))$ our  results agree with those of \cite
{13} obtained by an entirely different method, while for
$U_q(osp(2|2n))$ our results are new. In the limit $q\rightarrow 1$ our
results are also in agreement with \cite {14} taking into account the
different conventions adopted.

\subheading {Acknowledgement}

This research is supported by the Australian Research Council.
YZZ is financially supported by the Kyoto University Foundation.

\heading Appendix \endheading

Here we will prove proposition 2. To do so, we will require a few
technical results. The first is due to Kostant \cite {22}.
\proclaim{Lemma A1} Impose a partial ordering on $H^*$ with the order
relation $\mu\geq\nu$ if $\mu -\nu\in D^+$. Let $\nu_0\in D^+$ be
arbitrary and let $f$ be any polynomial function on $H^*$. If $f(\nu)=0$
for all $\nu \geq \nu_0$ then $f$ vanishes identically on $H^*$.
\endproclaim

Lemma A1 is proved in \cite {22} for the case of non $\Bbb Z_2$-graded
algebras, but is easily extended to the superalgebra case. We will also
require the following \cite {21}.
\proclaim{Lemma A2} Let $U_q(g_0)\subset U_q(g)$ be the (non-graded)
quantum algebra generated by $\{e_i,\,f_i,\,h_i\}_{i=1}^l$. If
$$V(\La)\otimes V(\mu)=\bigoplus_im_i V(\nu_i),$$
where $\mu,\,\nu_i$ are all typical and dominant then we have the $U_q(g_0)$
decomposition
$$V(\La)\otimes V_0(\mu)=\bigoplus_im_i V_0(\nu_i),$$
where $V_0(\nu)$ denotes the $U_q(g_0)$ module of highest weight $\nu$ and
$V(\La)$ is a direct sum of such modules; viz.
$$V(\La)=\bigoplus_jV_0(\La_j). $$
\endproclaim

In the case that $\La$ is subordinate to $\mu$ it follows from lemma A2
that each $\La_j$ is subordinate to $\mu$ when considered as
$U_q(g_0)$ weights. Since in the non-graded case for
$\mu\in D^+_{\La_j},\,\forall j$,
the multiplicity of $V_0(\mu +\la_i)\subset V(\La)\ot V_0(\mu)$ is equal
to the multiplicity of $\la_i\in\Pi(\La)$ \cite {22}, it also follows
from lemma A2 that the multiplicity is the same in the $\Bbb Z_2$-graded
decomposition (8).

We now recall the
following result due to Parthasary, Ranga-Rao and Varadarajan \cite {23}.
\proclaim{Lemma A3} For $\la_i\in\Pi(\La)$ the multiplicity of
$V_0(\mu+\la_i)\subset V(\La)\ot V_0(\mu)$ is given by $dim\,\,
V_{\la_i,\mu}(\La)$ where
$$V_{\la_i,\mu}(\La)=\{v\in V_{\la_i}(\La)|e_j^{ <\mu +\rho , \al_j>}v=0,
\,1\leq j\leq l\}.$$
Above, $V_{\la_i}(\La)\subset V(\La)$ is the space of all vectors of
weight $\la_i$.  \endproclaim

It is apparent from lemma A3 that $\La$ is subordinate to $\mu$ if and
only if
$$\align
e_i^{<\mu+\rho,\al_i>}v=&e_i^{<\mu,\al_i>+1}v   \\
                       =&0,\qquad   \forall 1\leq i\leq l,\quad v\in
		       V(\La).\endalign$$
In view of the finite dimensionality of $V(\La)$ there must exist
positive integers $m_1,.....,m_l$ such that
$$e_i^{m_i+1}v=0,\qquad\forall\,v\in V(\La). \tag {10} $$
Now let $\om_1,......,\om_l$ denote the fundamental weights of $U_q(g)$.
For $i\neq 0$ they are defined by
$$<\om_i,\,\al_j>=\delta_{ij},$$
while $\om_0$ is defined by
$$(\om_0,\,\al_i)=\delta_{i0}.$$
Setting
$$\nu_0=\sum_{i=0}^lm_i\om_i$$
where $m_i,\, 1\leq i\leq l$ are to satisfy (10) and $m_0$ is arbitrary,
it is clear that $\La$ is subordinate to $\nu_0$. If $\nu\in D^+$
satisfies $\nu \geq \nu_0$ then $\La$ is also subordinate to $\nu_0$.
We have thus shown
\proclaim {Lemma A4} Let $\La\in D^+.$ There exists $\nu_0\in D^+$ such
that $\La$ is subordinate to $\nu$ whenever $\nu \geq \nu_0$.
\endproclaim

To see that proposition 2 follows from lemmas A1-4, it suffices to
observe that the set of typical $\mu\in D^+$ are Zariski dense in $D^+$.
(For a description of the Zariski topolgy on $H^*$ see \cite {24}.)
Since polynomial functions are continuous in the Zariski topology, we
may conclude that proposition 2 holds.

\+\cr
\parindent 8mm

\heading References \endheading

\ref\no 1\by M. Jimbo\jour Lett. Math. Phys. \vol 10 \yr 1985 \page 63
\moreref \jour ibid. \vol 11 \yr 1986 \page 247 \moreref \jour V.G.
Drinfeld, Proc. I.C.M. \yr Berkeley, 1986 \page 798 \endref

\ref \no 2\by P.P. Kulish and N. Yu Reshetikhin\jour Lett. Math. Phys.
\vol 18 \yr 1989 \pages 143
\moreref\jour  M. Chaichian and P.P. Kulish, Phys. Lett. \vol B236 \yr
1990\pages 242 \moreref \jour A.J. Bracken, M.D. Gould and R.B. Zhang,
Mod. Phys. Lett. \vol A5 \yr 1990 \pages 331
\moreref \jour R. Floreanini, V.P. Spiridonov and L. Vinet,
Commun. Math. Phys. \vol 137 \yr 1991 \pages 149
\moreref \jour R. Floreanini, D.A. Leites and L. Vinet, Lett. Math.
Phys.\vol 23 \yr 1991 \pages 127 \moreref \jour M. Scheunert, Lett.
Math. Phys. \vol 24 \yr 1992 \pages 173 \moreref \jour M. Scheunert,
J. Math. Phys. \vol 34 \yr 1993 \page 3780\endref

\ref \no 3 \by S.M. Khoroshkin and V.N. Tolstoy \jour Commun. Math.
Phys. \vol 144 \yr 1991 \pages 599\endref

\ref\no 4\by H. Yamane \jour Proc. Jpn. Acad. \vol A67 \yr 1991 \page
108 \endref

\ref\no 5\by S. Sarker\jour J. Phys. A: Math. Gen \vol 23 \yr 1990 \page
L409\moreref\jour ibid. \vol 24 \yr 1991 \page 1137 \moreref \jour
F.H.L. Essler and V.E. Korepin, Phys. Rev. \vol B46 \yr 1992 \page 9147
\endref

 \ref\no 6\by S.-M. Fei and R.-H. Yue \jour J. Phys. A: Math. Gen. \vol
 27 \yr 1994 \page 3715 \moreref \jour A. Foerster and M. Karowski,
 Nucl. Phys. \vol B408 \yr 1993 \page 512 \endref

 \ref\no 7\by F.H.L. Essler, V.E. Korepin and K. Schoutens \jour Phys.
 Rev. Lett. \vol 68 \yr 1992 \page 2960\moreref\jour ibid. \vol 70 \yr
 1993 \page 73\endref

 \ref\no 8\by A.J. Bracken, M.D. Gould, J.R. Links and Y.-Z. Zhang
 \jour Phys. Rev. Lett. \vol 74 \yr 1995 \page 2768 \endref

 \ref\no 9\by L.D. Faddeev, N. Yu Reshetikhin and L.A. Takhtajan \jour
 Algebraic Analysis \vol 1 \yr 1989 \page 129 \moreref \jour M.D. Gould,
 R.B. Zhang and A.J. Bracken, J. Math. Phys. \vol 32 \yr 1991 \page 2298
 \moreref \jour J.R. Links and M.D. Gould, Rep. Math. Phys. \yr 1992
 \page 91 \vol 31 \endref

\ref \no 10\by R.B. Zhang and M.D. Gould \jour J. Math. Phys.\vol 32 \yr
1991 \pages 3261 \endref

\ref\no 11 \by M.L. Nazarov \jour Lett. Math. Phys. \vol 21 \yr 1991
\page 123\moreref\jour R.B. Zhang, {\it The quantum super Yangian and
Casimir operators of} $U_q(gl(m|n))$, Lett. Math. Phys. \toappear
\endref

\ref\no 12 \by R.B. Zhang\jour J. Math. Phys. \vol 33 \yr 1992 \page
3918\moreref \jour J.R. Links, M.D. Gould and R.B. Zhang, Rev. Math.
Phys. \vol 5\yr 1993 \page 345\moreref\jour M.D. Gould, I. Tsohantjis
and A.J. Bracken, Rev. Math. Phys. \vol 5 \yr 1993 \page 533 \endref

\ref \no 13 \by J. R. Links and R.B. Zhang \jour J. Math. Phys. \vol
34\yr 1993 \pages 6016\endref

\ref \no 14 \by A.M. Bincer \jour J. Math. Phys. \vol 24\yr 1983 \pages
2546 \moreref
\jour  M. Scheunert,  J. Math. Phys. \vol 24 \yr 1983 \pages
2681 \moreref\jour H.S. Green and P.D. Jarvis, J. Math. Phys. \vol 24
\page 1681 \yr 1983\endref

\ref\no 15 \by M.D. Gould and M. Scheunert\jour J. Math. Phys. \vol 36
\page 435 \yr 1995\moreref \jour J.R. Links and M.D. Gould, J. Math.
Phys. \vol 36 \yr 1995 \page 531 \endref

\ref\no 16 \by A. J. Bracken, G.W. Delius, M.D. Gould and Y.-Z Zhang
\jour J. Phys. A: Math. Gen \vol 27 \yr 1994 \page 6551 \moreref\jour
G.W. Delius, M.D. Gould, J.R. Links and Y.-Z. Zhang, {\it On type I
quantum affine superalgebras}, hep-th/9408006, Int. J. Mod. Phys.
A \toappear  \moreref \jour G.W. Delius, M.D.
Gould, J.R. Links and Y.-Z. Zhang, {\it Solutions of the Yang-Baxter
equation with extra non-additive parameters II:} $U_q(gl(m|n))$,
hep-th/9411241 \endref

\ref\no 17\by J.R. Links and M.D. Gould \jour Lett. Math. Phys. \vol 26
\yr 1992 \page 187\moreref \jour J.R. Links, M. Scheunert and M.D.
Gould, Lett. Math. Phys. \vol 32 \yr 1994 \page 231 \endref

\ref\no 18\by M.D. Gould, J.R. Links and Y.-Z. Zhang \paper Type I
quantum superalgebras, $q$-supertrace and two-variable link
polynomials\jour University of Queensland and Kyoto University
preprint \yr 1995 \endref

\ref\no 19\by V.G. Kac\jour Lect. Notes in Math. \vol 676 \yr 1978 \page
597 \endref

\ref\no 20\by R.B. Zhang\jour J. Math. Phys. \vol 34 \yr 1993 \page 1236
\moreref \jour J. Phys. A: Math. Gen \vol 26 \yr 1993 \page 7041 \endref

\ref\no 21\by M.D. Gould and J.R. Links \jour in preparation \endref

\ref\no 22\by B. Kostant\jour J. Funct. Anal. \vol 20 \yr 1975 \page 257
\endref

\ref\no 23\by K.R. Parthasary, R. Ranga-Rao, and V.S. Varadarajan\jour
Ann. Math. \vol 85 \yr 1967 \page 383 \endref

\ref\no 24\by J.E. Humphreys \paper Introduction to Lie algebras and
representation theory\jour Springer-Verlag, Berlin \yr 1972 \endref

\enddocument